\documentclass{article}
\usepackage{amssymb}
\usepackage{graphicx}

\usepackage{subcaption}
\usepackage{amssymb,latexsym,amsmath,amsthm}
\usepackage{xcolor}

\usepackage{tikz,tikz-cd}
\usetikzlibrary{hobby}
\usetikzlibrary{arrows}
\usetikzlibrary{positioning}
\usetikzlibrary{shapes,snakes}
\usetikzlibrary{cd}

\def\Rnn{\mathbb{R}^{n\times n}}
\def\R{\mathbb{R}}
\def\bone{{\bf 1}}

\def\yvec{{\bf y}}
\def\xsym{\boldsymbol{x}}
\def\ysym{\boldsymbol{y}}
\DeclareMathOperator{\diag}{diag}

\newtheorem{theorem}{Theorem}[section]

\newtheorem{prop}{Proposition}[section]

\theoremstyle{definition}
\newtheorem{definition}[theorem]{Definition}
\newtheorem{example}[theorem]{Example}

\theoremstyle{remark}
\newtheorem{remark}[theorem]{Remark}

\numberwithin{equation}{section}

\begin{document}

\title{Dynamic Katz and Related Network Measures}

\author{
Francesca Arrigo\thanks{Department of Mathematics and Statistics, University of Strathclyde, Glasgow, UK, G1 1XH (\textit{francesca.arrigo@strath.ac.uk}).},
Desmond J. Higham\thanks{School of Mathematics, University of Edinburgh, James Clerk Maxwell Building, Edinburgh, UK, EH9 3FD (\textit{d.j.higham@ed.ac.uk}).},
Vanni Noferini\thanks{Aalto University, Department of Mathematics and Systems Analysis, P.O. Box 11100, FI-00076, Aalto, Finland (\textit{vanni.noferini@aalto.fi}).}
and 
Ryan Wood\thanks{Corresponding author. Aalto University, Department of Mathematics and Systems Analysis, P.O. Box 11100, FI-00076, Aalto, Finland (\textit{ryan.wood@aalto.fi}).}}

\date{}

\maketitle
\begin{abstract} 
    We study walk-based centrality measures for time-ordered network sequences. For the case of standard dynamic walk-counting, we 
    show how to derive and compute centrality measures induced by analytic functions. We also prove that dynamic Katz centrality, based on the resolvent function, has the unique advantage of 
    allowing computations to be performed entirely at the node level.
    We then consider two distinct types of backtracking and 
    develop a framework for capturing dynamic walk combinatorics 
    when either or both is disallowed. 
\end{abstract}


\section{Introduction}\label{sec:intro}
Centrality measures play an important role in network analysis 
by quantifying the importance of each node \cite{BV14,Estradabook}. 
How to make the notion of importance mathematically precise is highly context dependent, and hence a range of centrality measures exist. 
Walk-based centralities constitute 
a large and widely-used class, based on
counting the 
number of walks beginning from each node, with some form of
downweighting based on the walk length \cite{BEK13,EHSiamRev}. Different types of weights lead to different families of centralities, such as Katz versus exponential; one or more parameters allow for even further freedom within each family. Moreover, in recent years it has been proposed that one might wish to ignore walks that immediately backtrack, to obtain centrality measures that are better suited for certain applications. Such centrality measures 
are said to be \emph{non-backtracking} and are known to offer tangible benefits \cite{AGHN17a,AGHN17b,Beyond,GHN18}. In this paper, we present a
framework for analyzing and computing walk-based
centrality measures, including 
non-backtracking versions, in the case of  
 time-evolving networks.
 
 Key contributions in this work are as follows.
 \begin{itemize}
     \item A proof that Katz (or resolvent-based) centrality is unique among
     matrix function-based versions, in that a combinatorially 
     correct expression exists in the form of products of functions of node-level adjacency matrices (Theorem~\ref{thm:katzunique});
     thereby showing that dynamic Katz has a major computational 
      advantage.
     \item Edge-level expressions for other matrix function-based versions ((\ref{eq:f1}) and (\ref{eq:f2})) using a block matrix construction, with an accompanying convergence theory
     (Proposition~\ref{prop:generalf}).
     \item Corresponding computable expressions for non-backtracking versions of matrix function-based temporal centrality measures
      (
      \eqref{eq:nbtf1} and \eqref{eq:nbtf2}) along with an accompanying convergence theory 
      (Propositions~\ref{prop:radius} and \ref{prop:fnbt}).
 \end{itemize}

The material is organized as follows.
Section~\ref{sec:background}
sets up some notation and explains what we mean by a temporal
network. 
In section~\ref{sec:expono} we 
state and prove a unique property of dynamic Katz.
Sections~\ref{sec:dyncom} and \ref{sec:measures_temporal}
develop a framework in which 
general matrix function-based centralities can be handled.
In section~\ref{sec:nbt} we 
introduce and study new non-backtracking versions of these measures,
focusing on three main variations where backtracking is disallowed 
within timesteps, across timesteps, or both.
Section~\ref{sec:conc} concludes with a brief discussion.

\section{Background}\label{sec:background}
In this section we review some definitions and notation associated with graphs and 
some of their matrix representations.

A {\it network}, or {\it graph}, is an ordered pair of sets $G = (V,E)$, 
where $V$ is the set of {\it nodes}  
and $E\subset V\times V$ is the set of \emph{edges} among the nodes \cite{Estradabook,Newmanbook}. 
We consider graphs with 
unweighted edges and no self-loops.
Given an edge $i\to j$ we will call $i$ its {\it source} and $j$ its {\it target}. If the edge from 
$j\to i$ also exists, we will refer to this as the {\it reciprocal} of the edge $i\to j$ and we will say 
that edge $i\to j$ is {\it reciprocated}.

Throughout this work $I$ denotes the identity matrix, $\bone$ denotes the vector of all ones; the size of both, unless otherwise stated, is chosen such that the resulting expressions in which they appear are coherent. The symbol $\star$ will be used as a placeholder to denote an arbitrary vertex, where appropriate.

Two definitions that are central to our work those of 
a {\it walk} and a {\it (non-)backtracking walk} around a graph. 
We review them in the following. 
\begin{definition}
A {\rm walk of length $r$} is a sequence of $r$ 
edges $(e_1,e_2,\ldots,e_r)$ such that the target of $e_\ell$ coincides with 
the source of $e_{\ell+1}$ for all $\ell=1,2,\ldots, r-1$.
\end{definition}
\begin{remark}
An equivalent definition can be given in terms of the nodes: namely, a walk of length $r$ can also be seen as a sequence of $r+1$ nodes 
$(i_1,i_2,\ldots,i_{r+1})$ such that $(i_\ell,i_{\ell+1})\in E$ for all $\ell = 1,2,\ldots, r$. 
\end{remark}

\begin{definition}\label{def:nbt}
Let $G$ be a graph. A walk in $G$ is said to be {\rm backtracking} if any 
two consecutive edges in it are reciprocated. The walk is said to be {\rm non-backtracking (NBT)} otherwise.
\end{definition}

\subsection{Graphs and matrices}\label{ssec:matrices}
A widely used linear algebraic representation of a graph is its adjacency matrix.
\begin{definition}\label{def:A}
Let $G = (V,E)$ be an unweighted graph with $n$ nodes. Its {\rm adjacency matrix} $A \in \mathbb{R}^{n \times n}$ is entrywise 
defined as:
$$ A_{ij}=\begin{cases}
1 \ \ \text{if}\,\, (i,j) \in E\\
0 \ \ \text{otherwise}
\end{cases}
$$
for all $ i,j = 1,2,\ldots, n$.
\end{definition}

Other matrices of interest for the purpose of our work are the {\it source} and {\it target} (or {\it terminal}) 
matrices associated with $G$ and the adjacency matrix of the line graph of $G$; see, e.g., \cite{Estradabook,VDF09,ZV08}.

\begin{definition}\label{def:LR}
Let $G$ be an unweighted graph with $n$ nodes and $m$ edges. Its {\rm source} and {\rm target} (or {\rm terminal})  {\rm matrices} $L,R \in \mathbb{R}^{m \times n}$ are entrywise defined as:
$$L_{ei}=\begin{cases}
1 \ \ \text{if edge~} e \text{~has the form~} i\to\star\\
0 \ \ \text{otherwise}
\end{cases}
$$
and
$$
R_{ej}=\begin{cases}
1 \ \ \text{if edge~} e \text{~has the form~} \star\to j\\
0 \ \ \text{otherwise}
\end{cases}$$
respectively, for all $e = 1,2,\ldots,m$ and for all $ i,j = 1,2,\ldots, n$.
\end{definition}

The {\it line graph} (or {\it interchange graph} or {\it dual graph}) of $G$ is constructed from $G$ as follows: edges in the original graph are regarded as nodes 
and two nodes
$i \to j$ and $k \to \ell$ 
 in this new line graph 
are connected if $ j = k$, that is, if, together, they represent 
a walk of length two in the original graph $G$; see, e.g.,~\cite{ore62}.
The adjacency matrix $W\in\R^{m\times m}$ of the line graph can then be entrywise 
defined in terms of elements of $G$ as 
\[
W_{i\to j,k\to \ell} = \delta_{jk}, 
\]
where $\delta_{jk}$ is the Kronecker delta.

\subsubsection{Walk-based centrality measures}

Katz centrality \cite{Katz53} is a widely used centrality measure that assigns to each node $i$ the $i$th element of the vector
\begin{equation}\label{eq:Katz}
    \yvec(\alpha) = \sum_{r=0}^\infty \alpha^rA^r\bone.
\end{equation}
Combinatorially, Katz centrality can be understood as assigning to node $i$ the sum over all possible walk lengths of the number of walks originating from $i$ of length $r$ scaled by $\alpha^r$, where $\alpha>0$ is a downweighting parameter. 
Indeed, it is easy to show that $(A^r)_{ij}$ is the number of walks of length $r$ originating from node $i$ and ending at node $j$ and hence the above interpretation of the entries of $A^r\bone$ readily follows. 

Furthermore, whenever $0 < \alpha < 1/\rho(A)$, where $\rho(A)$ denotes the spectral radius of $A$, the series appearing in~\eqref{eq:Katz} converges and the vector of Katz centralities can be equivalently defined as 
\[
    \yvec(\alpha) = (I-\alpha A)^{-1}\bone.
\]

The idea underlying Katz centrality can be generalized to obtain other centrality measures for nodes defined in terms of matrix functions and their entries, or sums thereof; see, e.g.,~\cite{EHSiamRev} and references therein. 
Suppose that the analytic function 
\begin{equation}
f(z) = \sum_{r=0}^\infty c_r z^r
\label{eq:fz}
\end{equation}
has nonnegative Maclaurin coefficients $c_r \geq 0$ and is convergent for $|z|<R$. 
Then, for any $0<\alpha<R/\rho(A)$ the matrix power series
\[ f(\alpha A) = \sum_{r=0}^\infty c_r \alpha^r A^r\]
also converges~\cite{HigFM}. 
Entries and sums of entries of $f(\alpha A)$ can thus be interpreted in terms of the number of walks taking 
place in the network, with shorter walks being given more weight since $c_r \alpha^r\to 0$ as $r\to \infty$. 
Because of this combinatorial interpretation, the function $f$ may be used to define centrality measures of nodes:
\begin{itemize} 
\item {\it $f$-total communicability:} $f(\alpha A)\bone$; and 
\item {\it $f$-subgraph centrality}, defined as the diagonal of $f(\alpha A)$.
\end{itemize}

A notable special case is that of $c_r=1/r!$ or, equivalently, $f(z)=e^z$, which yields the 
popular {\it subgraph centrality} and {\it total communicability} measures for nodes, defined in terms of 
the diagonal elements and row sums of $e^A$, respectively~\cite{BK13,EV05}.

\subsection{Temporal Networks}\label{ssec:temporal_net}
The definition of a network may be extended in various ways to a time-evolving setting \cite{holme11}.
We consider here the case where, for a fixed collection of nodes, 
a network is recorded at a discrete set of time points.
This is relevant in many message passing and digital
interaction 
contexts; see, e.g., \cite{bassett11,SNAM2013,MH12}.
 
\begin{definition}
A {\it temporal network} $\mathcal{G}$ with $n$ nodes is a finite 
ordered sequence of networks $\mathcal{G} = (G^{[1]}, \dots, G^{[N]})$ associated with a finite increasing collection of time stamps ${\mathcal{T}}_N = (t_1, \dots, t_N)$, where each network $G^{[\tau]} = (V^{[\tau]},E^{[\tau]})$ exists at time $t = t_\tau$ for $\tau =1,2,\ldots,N$. 
\end{definition}

In the following we will assume that $V^{[\tau]}=\mathcal{V} = \{1,\ldots, n\}$ for all $\tau = 1,2,\ldots,N$; that is, the node set is fixed over time. In principle, nodes which join or leave 
the temporal network can be accommodated in this framework---they will be isolated nodes when inactive. 
We further assume that, for all $\tau = 1,2,\ldots, N$, $G^{[\tau]}$ is unweighted and contains no self-loops. We denote by $m_\tau=|E^{[\tau]}|$ the number of edges appearing in $G^{[\tau]}$. Further, we denote by $m=\sum_{\tau = 1}^N m_\tau$ the total number of edges 
in ${\mathcal{G}}$.

Each of the graphs $G^{[1]},G^{[2]},\ldots,G^{[N]}$ appearing in the definition of ${\mathcal{G}}$ can be described via their adjacency matrices, which we denote by $A^{[\tau]}$ for all $\tau = 1,2,\ldots,N$. 
These are implicitly used in~\cite{holme11} to define the {\it adjacency index} (or {\it presence function}~\cite{casteigts2012time}) of a temporal network ${\mathcal{G}}$, which is a third-order tensor ${\mathcal{A}}\in\R^{n\times n\times N}$ entrywise defined as
\[
{\mathcal{A}}_{ij\tau} = (A^{[\tau]})_{ij}
\]
for all $i,j\in {\mathcal{V}}$ and $\tau = 1,2,\ldots,N$.
We will instead work with the separate
adjacency matrices $A^{[\tau]}$;
furthermore, we will denote by $L^{[\tau]}$ and $R^{[\tau]}$
the corresponding source and target matrices, respectively, and by $W^{[\tau]}$ the adjacency matrix of the line graph of $G^{[\tau]}$.

The definition of a walk across a network can be extended to the setting of temporal graphs 
as follows. 
\begin{definition}\label{def:Twalk1}
{\rm A walk of length $r$} across a temporal network is defined as an ordered sequence of $r$ edges $(e_1, e_2, \dots, e_{r})$ such that the target of $e_\ell$ coincides with the source of $e_{\ell+1}$ for all $\ell =1,\ldots,r-1$ and, moreover, that $e_\ell \in E^{[\tau_1]}, e_{\ell+1} \in E^{[\tau_2]}$ for some $1\leq \tau_1 \leq \tau_2\leq N$.
\end{definition}
\begin{remark}
Equivalently, we could see a
 walk of length $r$ as an ordered sequence of $r+1$ nodes $(i_1, i_2, \dots, i_{r+1})$ such that for all $\ell = 2,\ldots,r$ it holds that $i_{\ell-1}\to i_{\ell} \in E^{[\tau_1]}$ and $i_\ell\to i_{\ell+1} \in E^{[\tau_2]}$ for some $1\leq\tau_1 \leq \tau_2\leq N$.
 \end{remark}

We stress that multiple edges can be used within a time stamp and, moreover, a walk is allowed to remain inactive   
for some of the time stamps. 
We will sometimes use the notation 
$i \xrightarrow[]{t_\tau} j$
to denote the edge $i\to j\in E^{[\tau]}$.

\bigskip

In the next section, after reviewing how Katz centrality was extended in~\cite{GHPE11} to the temporal setting 
by simply considering multiplications of matrix functions of the form $(I-\alpha A^{[\tau]})^{-1}$, 
we prove that Katz centrality is indeed the only centrality measure for which this approach 
can be employed. Generalizations to the temporal setting of other centrality measures, such as, e.g., 
total communicability, may be computed using the matrix construction that we present in 
section~\ref{sec:dyncom} in order to respect the combinatorics of walks. 

\section{Walk-based centrality measures in temporal networks}\label{sec:expono}
A generalization of Katz centrality to temporal networks is proposed in \cite{GHPE11} using the dynamic communicability matrix 
\begin{equation}\label{eq:dcm}
\mathcal{Q}  := \prod_{\tau=1}^N (I-\alpha A^{[\tau]})^{-1}.
\end{equation}  
This assigns to node $i$ the $i$th element of the vector
\begin{equation}
\ysym(\alpha) =\mathcal{Q}\bone=  (I-\alpha A^{[1]})^{-1}\cdots (I-\alpha A^{[N]})^{-1}\bone.
\label{eq:Qx}
\end{equation}
It can be checked that for $0<\alpha<(\max_\tau \{\rho(A^{[\tau]})\})^{-1}$, 
the component 
$\ysym(\alpha)_i$ corresponds to the weighed sum of all walks 
that emerge from node $i$, where walks of length $r$ are weighted as $\alpha^r$; see~\cite{GHPE11} for more details.

We note that this combinatorial interpretation, where walks 
of length $r$ are weighted by the appropriate coefficient $c_r$ in \eqref{eq:fz}, may fail 
if we were to consider matrix functions other than the resolvent in 
the expression \eqref{eq:Qx}. 
For example, 
the exponential case 
$c_r = \beta^r/r!$ in 
\eqref{eq:fz},
for some $\beta >0$,
was extended to the time-dependent setting 
using a quantum physics motivation 
to give \cite{estrada2013temporal}
\[
\prod_{\tau = 1}^Ne^{\beta A^{[\tau]}} = e^{\beta A^{[1]}}e^{\beta A^{[2]}}\cdots e^{\beta A^{[N]}}.
\]
However, \emph{from a combinatorial viewpoint}
this generalization has the drawback of not weighting walks 
consistently with their length. Consider for example a walk of length five from $i$ to $j$ that is realized by traveling two edges at time $t=t_1$, one edge at time $t=t_2$, and two edges at time $t=t_3$. 
Because of its length, we would expect this walk to be weighted as $\frac{\beta^5}{5!}$ in keeping with the theory that we have seen for static networks. 
However, it is straightforward to check that the $(i,j)$ entry of $e^{\beta A^{[1]}} e^{\beta A^{[2]}}e^{\beta A^{[3]}}$ is weighted as $\frac{\beta^2}{2!}\frac{\beta^1}{1!}\frac{\beta^2}{2!}=\frac{\beta^5}{4}$. 
Furthermore, this walk of length five is weighted differently from another walk of length 
five that travels, say, 
one edge at times $t=t_1, t_2$ and three edges at time $t=t_3$. Indeed, the weight of 
this second walk is $\frac{\beta}{1!}\frac{\beta}{1!}\frac{\beta^3}{3!}=\frac{\beta^5}{6}$. 
Thus, a simple generalization to the temporal setting as the one presented in \cite{estrada2013temporal} 
not only fails to respect the combinatorics of walks, with walks of length $r$ not weighted as 
$\frac{\beta^r}{r!}$, but it may also be inconsistent 
in weighting different walks of the same length. 

\smallskip

In the following subsection we show that resolvent-based 
centrality measures are unique in this regard: they are the only 
functions that respect the combinatorics of walks when translated to the 
temporal setting with the simple expression
\[
f(\alpha A^{[1]})f(\alpha A^{[2]})\cdots f(\alpha A^{[N]}).
\]

\subsection{A remarkable property of Katz centrality}\label{ssec:Katz_unique}

Is it possible for choices of $f$ other than the resolvent (including the popular $f(z)=e^z$) to compute $f$-centralities on a time-evolving network directly using products of functions of adjacency matrices? 

While the example above clarified that the exponential subgraph and total communicability on a time-evolving network cannot be computed by simply multiplying the matrix exponential of each adjacency matrix, it is still {\it a priori} possible that they could be computed by multiplying some \emph{other} function of each adjacency matrix, say, via $\prod_{\tau=1}^N g(\alpha A^{[\tau]})$ for some $g(z)$. More generally, posing this question for $f$-centrality measures is equivalent to the following combinatorial problem: Given the function $f(z) = \sum\limits_{r=0}^\infty c_r z^r$, $c_r \geq 0$, find a function $g(z)$ such that
 \begin{equation}\label{eq:Estradaiswrong}
\sum_{r=0}^\infty c_r \alpha^r h_r(x_1,\dots,x_N)=  \prod_{i=1}^N g(\alpha x_i),
\end{equation} 
where 
\[ h_r(x_1,\dots,x_N)= \sum_{i_1+\dots+i_N=r} x_1^{i_1} \dots x_N^{i_N}\]
is the $r$th complete homogeneous symmetric polynomial in $N\geq 2$ variables. Indeed, finding a solution to this scalar problem would imply that the weights assigned to a walk in the matrix $\prod_{\tau=1}^N g(\alpha A^{[\tau]})$ 
\begin{itemize}
    \item depend only on the total walk length and not, for example, on how it is split between time frames; and
    \item hold for any possible time evolving graph.
\end{itemize}

It turns out that this is essentially (up to a multiplicative constant and a linear change of variable) only possible if $f(z)$ is the resolvent. Since multiplication by a constant does not change the ranking while a linear change of variable is tantamount to changing the Katz parameter $\alpha$, it follows that  no $f$-centrality other than Katz (i.e., resolvent subgraph) centrality is applicable to time-evolving networks via an expression only based on products of functions of the adjacency matrices. For example, there is no way to correctly compute the combinatorics of exponential centrality on time-evolving graphs by directly multiplying some functions of the adjacency matrix.\footnote{Of course, this does not mean that computing exponential centralities on time evolving graphs is not possible: indeed, using the multilayer approach proposed in this paper is for example a combinatorially exact, albeit potentially demanding in terms of computational complexity, way to compute it.} The theorem below gives a formal proof of this claim.

\begin{theorem}\label{thm:katzunique}
Given the sequence of nonnegative real numbers $(c_r)_r$, let $f(z)=\sum_{r=0}^\infty c_r z^r$ and let $N \geq 2$. The functional equation \eqref{eq:Estradaiswrong} has a solution $g(z)$ if and only if there exist nonnegative constants $\gamma,\delta \geq 0$ such that $c_r=\gamma \delta^r$ for all $r$, i.e, if and only if $f(z)=\gamma(1-\delta z)^{-1}$. Moreover, in that case the solution is \[ g(z)= \sqrt[N]{\gamma}\ (1-\delta z)^{-1}. \]
\end{theorem}

\begin{proof}
Suppose first that  $c_r=\gamma\delta^r$ for some $\gamma,\delta\geq 0$, then  (see also \cite{GHPE11} and \cite[Ch. 7]{Stanley})
\[
\sum_{r=0}^\infty c_r \alpha^r h_r(x_1,\dots,x_N)  = \gamma  \sum_{r=0}^\infty \alpha^r h_r(\delta x_1,\dots,\delta x_N)  = 
\gamma \prod_{i=1}^N \sum_{r=0}^\infty (\delta x_i \alpha)^r = \gamma \prod_{i=1}^N \frac{1}{1-\alpha \delta x_i},
\]
so taking
$ g(\alpha x_i) = \sqrt[N]{\gamma} (1-x_i \delta \alpha)^{-1}$
yields a solution to \eqref{eq:Estradaiswrong}.

Conversely suppose that \eqref{eq:Estradaiswrong} has a solution $g(z)$, and assume that $f(0)=1$. In doing so we do not lose  generality, for it is clear that $g(z)$ solves \eqref{eq:Estradaiswrong} for $f(z)$ if and only if, for all $\gamma \geq 0$, $\sqrt[N]{\gamma}\ g(z)$ solves it for $\gamma f(z)$.  Recall now the following (classic) property of complete homogeneous symmetric polynomials \cite[Sec. 7.5]{Stanley}:
\[ h_r(\underbrace{1,\dots,1}_{\ell \text{ variables}},\underbrace{0,\dots,0}_{(N-\ell) \text{ variables}}) = \#\{ \text{monomials of degree } r \text{ in } \ell \text{ variables } \} = {\ell+r-1 \choose r}\]
Then, evaluating \eqref{eq:Estradaiswrong} at $(x_1,	\dots,x_N)=(0,\dots,0)$ yields $g(0)=1$. 
Next, we evaluate \eqref{eq:Estradaiswrong} at $(x_1,x_2,\dots,x_N)=(1,0,\dots,0)$ and conclude that $f(z)=g(z)$. Finally, by evaluating \eqref{eq:Estradaiswrong} at $(x_1,x_2,x_3,\dots,x_N)=(1,1,0,\dots,0)$, we obtain
\begin{equation}\label{eq:Eir} 
\sum_{r=0}^\infty (r+1) c_r \alpha^r = \sum_{r=0}^\infty \alpha^r \sum_{k=0}^r c_k c_{r-k}. 
\end{equation}
Now define $\delta:=c_1$ and proceed by induction.
We begin with $c_1=\delta^1$. Let us then assume that $c_k=\delta^k$ for all $1 \leq k \leq r-1$; equating the  coefficients of $\alpha^r$ in~\eqref{eq:Eir} we get
\[ (r+1) c_r = \sum_{k=0}^r c_k c_{r-k} = 2 c_r + \sum_{k=1}^{r-1}  \delta^r \quad\Longrightarrow\quad c_r = \delta^r.\]
\end{proof}

Motivated by this result, in sections~\ref{sec:dyncom} and \ref{sec:measures_temporal} we introduce a 
block upper-triangular matrix $M\in\R^{m\times m}$ which will allow us to 
generalize $f$-subgraph centrality and $f$-total communicability to 
treat temporal networks while respecting the combinatorics of walks. 
These results will then be generalized to the non-backtracking setting in section~\ref{sec:nbt}.

\section{The global temporal transition matrix}
\label{sec:dyncom}

We start by re-interpreting the results in \cite{GHPE11} via the line graph construction. 

To build intuition, suppose for simplicity that $N=2$ and let $A^{[1]}$ and $A^{[2]}$ be the adjacency matrices of the graphs appearing at time 
stamps $t_1 < t_2$ in the temporal network $\mathcal{G}$. 
For any adjacency matrix $A$ and adjacency matrix $W$ of the associated 
line graph, it holds that $A = L^T R$ and $W = RL^T$. 
These immediately yield

\begin{equation}\label{eq:Ar}
(A^r)_{ij} = (L^TW^{r-1}R)_{ij}
\end{equation} 
for all $i,j\in V$ and for all $r\geq 1$; see~\cite{Beyond,VDF09,ZV08}. 
\begin{remark}
Equation~\eqref{eq:Ar} above summarizes the fact that taking $r\geq 1$
steps in the node space, i.e., on the original 
graph $G$, corresponds to taking $r-1$ steps in the edge space, i.e., on the line graph associated with $G$. 
\end{remark}
Equation~\eqref{eq:Ar} together with~\eqref{eq:dcm} and the assumption that 
$0<\alpha<(\max_{\tau = 1,2}\{\rho(A^{[\tau]})\})^{-1}$, yields 
\begin{align*}
    \mathcal{Q} &= (I-\alpha A^{[1]})^{-1} (I-\alpha A^{[2]})^{-1}\\
    &= \left(I + \sum_{r = 1}^\infty \alpha^r {A^{[1]}}^r\right)\left(I + \sum_{r = 1}^\infty \alpha^r {A^{[2]}}^r\right)\\
    &=(I + \sum_{r = 1}^\infty \alpha^r {L^{[1]}}^T {W^{[1]}}^{r-1} R^{[1]})(I + \sum_{r = 1}^\infty \alpha^r {L^{[2]}}^T {W^{[2]}}^{r-1} R^{[2]})\\
    &= (I_n - \alpha {L^{[1]}}^T (I_m-\alpha W^{[1]})^{-1}R^{[1]})(I_n - \alpha {L^{[2]}}^T (I_m-\alpha W^{[2]})^{-1}R^{[2]}),
\end{align*}
and hence
\begin{equation}\label{eq:goodexample}
 \mathcal{Q}= I_n   + \alpha \left(  \sum_{\tau=1}^2 {L^{[\tau]}}^T (I_m - \alpha W^{[\tau]})^{-1} R^{[\tau]}\right) + 
\alpha^2 {L^{[1]}}^T (I_m-\alpha W^{[1]})^{-1} W^{[1,2]} (I_m- \alpha W^{[2]})^{-1} R^{[2]},
\end{equation}
where $W^{[1,2]}:=R^{[1]} (L^{[2]})^T$.
The combinatorial role of each term appearing in~\eqref{eq:goodexample} is as follows:
\begin{itemize}
\item $I_n$ counts walks that do not involve any edge within either graph, i.e.,  walks of length zero.
\item $\alpha \left(\sum_{\tau=1}^2 {L^{[\tau]}}^T (I - \alpha W^{[\tau]})^{-1} R^{[\tau]}\right)$ counts walks that use at least one 
edge within one of the two graphs, and no edges within the other: note that necessarily these walks have length at least one.
\item $\alpha^2 {L^{[1]}}^T (I_m - \alpha W^{[1]})^{-1} W^{[1,2]} (I_m- \alpha W^{[2]})^{-1} R^{[2]}$ counts walks that use at least one edge within each graph; note that necessarily these walks have length at least two, and, moreover, that at some point necessarily the walk has passed over an edge that exists at time $t_1$ to a subsequent edge that exists at time $t_2$, with the target of the time $t_1$ edge 
matching the source of the time $t_2$ edge (this is encoded in $W^{[1,2]}$).
\end{itemize}

Consider now the following block-triangular matrix:
\begin{equation}\label{eq:Mzero}
    M^{[1,2]} = \begin{bmatrix} W^{[1]} & W^{[1,2]} \\ 0 & W^{[2]}\end{bmatrix}.
\end{equation}
The following result holds. 
\begin{prop}\label{prop:properblocksW}
Suppose $M^{[1,2]}$ is defined as in~\eqref{eq:Mzero}. Then, when $0<\alpha<(\max_{\tau = 1,2}\{\rho(A^{[\tau]})\})^{-1}$ ,
  \[ (I - \alpha M^{[1,2]})^{-1} 
 = \begin{bmatrix} (I-\alpha W^{[1]})^{-1} & \alpha(I-\alpha W^{[1]})^{-1}W^{[1,2]}(I-\alpha W^{[2]})^{-1} \\ 0 & (I-\alpha W^{[2]})^{-1}\end{bmatrix}.\]
\end{prop}
\begin{proof}
Note that the $r$th power of $M^{[1,2]}$ is given by 
\[
 {M^{[1,2]}}^r = \begin{bmatrix} W^{[1]} & W^{[1,2]}\\
0 & W^{[2]}
\end{bmatrix}^r =  \begin{pmatrix}
{W^{[1]}}^r & \sum_{s = 0}^{r-1} {W^{[1]}}^s W^{[1,2]}{W^{[2]}}^{r-1-s} \\ 
0 & {W^{[2]}}^r\end{pmatrix}. \]
 Thus
\[
\sum_{r = 0}^\infty \alpha^r {M^{[1,2]}}^r = I_{2m}+ \sum_{r = 1}^\infty \begin{pmatrix}\alpha^r{W^{[1]}}^r & \alpha^r\sum_{s = 0}^{r-1} {W^{[1]}}^s W^{[1,2]}{W^{[2]}}^{r-1-s} \\ 0 & \alpha^r{W^{[2]}}^r \end{pmatrix}.
\]
Since the spectral radius of $W^{[\tau]}$ coincides with that of $A^{[\tau]}$ by Flanders' theorem, it follows that, when $0<\alpha<(\max_{\tau = 1,2}\{\rho(A^{[\tau]})\})^{-1}$, the diagonal blocks of the matrix on the right-hand side equal $(I-\alpha W^{[\tau]})^{-1}$, for $\tau = 1,2$, respectively. 

Finally, the non-zero off-diagonal block is such that
\begin{align*}
\sum_{r = 1}^\infty \alpha^r\sum_{s = 0}^{r-1} {W^{[1]}}^s W^{[1,2]}{W^{[2]}}^{r-1-s} &=  
\alpha\sum_{r = 0}^\infty \alpha^r\sum_{s = 0}^{r} {W^{[1]}}^s W^{[1,2]}{W^{[2]}}^{r-s} \\
&= \alpha \left(\sum_{s = 0}^\infty \alpha^s {W^{[1]}}^s \right) W^{[1,2]}\left(\sum_{s = 0}^\infty \alpha^s {W^{[2]}}^s \right)\\ 
&= (I_m-\alpha W^{[1]})^{-1}\alpha W^{[1,2]}(I_m-\alpha W^{[2]})^{-1}.
\end{align*}
The above, together with the fact that 
\[(I_{2m}-\alpha M^{[1,2]})^{-1} = I_{2m} + \sum_{r = 1}^\infty \alpha^r {M^{[1,2]}}^r,\]  
concludes the proof.
\end{proof}

If we now let 
\[\mathcal{L} = \begin{bmatrix}L^{[1]} \\ L^{[2]} \end{bmatrix}\quad \text{ and }\quad \mathcal{R} = \begin{bmatrix}R^{[1]} \\ R^{[2]} \end{bmatrix}
\]
it follows from 
(\ref{eq:Ar}) and Proposition~\ref{prop:properblocksW}
 that 
\[
{\mathcal{Q}} = I + \alpha {\mathcal{L}}^T (I-\alpha M^{[1,2]})^{-1}{\mathcal{R}}.
\]

This alternative description of ${\mathcal{Q}}$ is not necessary in practice, since ${\mathcal{Q}}$ can be 
easily computed by working in the node space, i.e., 
directly from the adjacency matrices $A^{[\tau]}$; however, 
this construction, appropriately generalized, will allow us to
define and compute 
 $f$-centrality measures for temporal networks for all functions $f$  
and NBT versions for temporal networks. 

These two points will be addressed in sections~\ref{sec:measures_temporal} and \ref{sec:nbt} below. 
In the reminder of this section, we will extend the definition of the matrices $M=M^{[1,2]}$, ${\mathcal{L}}$, 
and ${\mathcal{R}}$ to treat the case of $N \ge 2$ time stamps.  
To gain further intuition, we will also describe how the generalization of the matrix $M$ can be interpreted as the 
adjacency matrix of a static multilayer network.

\smallskip

We begin by extending the definitions and results above to the case of $N \ge 2$. 

\begin{definition}\label{def:Wt1t2}
Let ${\mathcal{G}} = (G^{[1]},G^{[2]},\ldots,G^{[N]})$ be a temporal network with $N$ time stamps $t_1 < t_2 < \cdots < t_N$, with  $G^{[\tau]} = ({\mathcal{V}}, E^{[\tau]})$.  
The {\rm backtrack-allowing transition matrix} $W^{[\tau_1,\tau_2]}$
where $1\leq\tau_1 < \tau_2\leq N$, is defined entrywise as follows:
\[
(W^{[\tau_1,\tau_2]})_{i \rightarrow j, k \rightarrow \ell} = \delta_{jk}
\]
for all $i\to j\in E^{[\tau_1]}$ and $k\to\ell\in E^{[\tau_2]}$. 
\end{definition}
These matrices are generally non-square and encode in their entries the presence 
of temporal walks of length two 
that take place across two distinct, but not necessarily consecutive, time stamps.

\begin{prop}\label{prop:Wt1t2}
Let ${\mathcal{G}} = (G^{[1]},G^{[2]},\ldots,G^{[N]})$ be a temporal network with $N$ time stamps $t_1< t_2<\cdots < t_N$. 
Let moreover $L^{[\tau]}$ and $R^{[\tau]}$ be the source and target matrices associated with graph $G^{[\tau]}$. 
Then for all $\tau_1 < \tau_2$ we have 
\[
W^{[\tau_1,\tau_2]} = R^{[\tau_1]}{L^{[\tau_2]}}^T.
\]
\end{prop}
\begin{proof}
The result may be established straightforwardly by working entrywise.
\end{proof}

We now extend the definition of $M^{[1,2]}$ to the case of $N\geq 2$ as follows. 
\begin{definition}\label{def:Mnolla}
Let ${\mathcal{G}} = (G^{[1]},G^{[2]},\ldots,G^{[N]})$ be a temporal network with $N$ time stamps 
$t_1< t_2<\cdots< t_N$.
The {\rm global temporal transition matrix} associated with ${\mathcal{G}}$ is the $m\times m$ block matrix:
\[
M:=M^{[1,2,\dots,N]} = \begin{bmatrix}  
W^{[1]} & W^{[1,2]} & W^{[1,3]} & \dots & W^{[1,N]} \\ 
0 & W^{[2]} & W^{[2,3]} & \dots & W^{[2,N]} \\
\vdots & & \ddots &  & \vdots \\
\vdots & & & \ddots & \vdots\\
0 & \dots & \dots & 0 & W^{[N]}
\end{bmatrix}.
\]
\end{definition}
We mention in passing that this is the adjacency matrix of a multilayer graph associated with the line graph of ${\mathcal{G}}$; see subsection~\ref{sec:multilayer} for more details.  
Similar constructions, although at node level, have been used over the years as a means to incorporate the temporal aspect in block-matrix 
representations of ${\mathcal{G}}$; 
examples include~\cite{al2021block,fenu2017block,MRMPO10,TMCPM17}. 
We note however that the adjacency matrix presented in Definition~\ref{def:Mnolla} does not correspond to the adjacency matrix of the line graph of any of the networks underlying the matrices presented in these references. Indeed, the number of nonzeros in the matrices presented in these references is strictly larger than the number of edges in the temporal graph, due to the presence of ``artificial edges" described in the off-diagonal blocks. Therefore the adjacency matrix of the associated line graph would be of larger dimension than the matrix in Definition~\ref{def:Mnolla}

 \begin{definition}\label{def:globalLR}
Let ${\mathcal{G}} = (G^{[1]},G^{[2]},\ldots,G^{[N]})$ be a temporal network with $N$ time stamps $t_1< t_2<\cdots< t_N$ 
and let $L^{[\tau]}$ and $R^{[\tau]}$ be the source and target matrices of $G^{[\tau]}$ for $\tau = 1,2,\ldots, N$. 
The {\rm global source} and {\rm global target} matrices to $N$ time stamps associated to ${\mathcal{G}}$ are
\[\mathcal{L}=\mathcal{L}^{[1,\dots,N]} = \begin{pmatrix} L^{[1]} \\ L^{[2]} \\ \vdots \\ L^{[N]}\end{pmatrix} \quad \text{ and } \quad \mathcal{R}=\mathcal{R}^{[1,\dots,N]} = \begin{pmatrix}R^{[1]} \\ R^{[2]} \\ \vdots \\ R^{[N]}\end{pmatrix},\]
respectively.
\end{definition}

It can be verified that, for all $\alpha\in(0,(\max_{\tau}\{\rho(A^{[\tau]})\})^{-1})$,
 Katz centrality 
 \eqref{eq:Qx} 
 on a time-evolving network with global temporal transition matrix $M=M^{[1,\dots, N]}$ can be rewritten as 
\begin{equation}\label{eq:ysym_N}
\ysym(\alpha) = 
{\mathcal{Q}}\bone =  [I + \alpha \mathcal{L}^T (I-\alpha M)^{-1} \mathcal{R}] \bone = \bone + \alpha \mathcal{L}^T (I-\alpha M)^{-1}\bone,    
\end{equation}
where we have used the fact that ${\mathcal{R}}\bone = \bone$. 

As we mentioned earlier in the section, this construction is not actually needed 
to compute $\ysym(\alpha) ={\mathcal{Q}}\bone$. However, it will be instrumental to 
the generalizations to $f$-centrality measures and their NBT equivalent. 
One observation that emerges from the above is that walks around the temporal graph 
${\mathcal{G}}$ can be enumerated by counting walks around the graph underlying $M$ and then 
projecting to the node space appropriately. 
We explore this connection in the next subsection. 

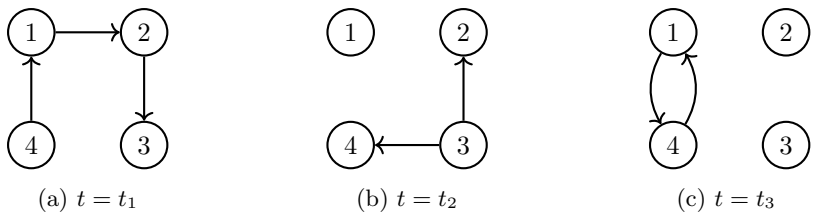
\begin{figure}
     \centering
     \begin{subfigure}[b]{0.3\textwidth}
         \centering
         \begin{tikzpicture}[]
    \node[circle,draw=black,thick] (11) at (0,1.5) {1};
    \node[circle,draw=black,thick] (21) at (1.5,1.5) {2};
    \node[circle,draw=black,thick] (31) at (1.5,0) {3};
    \node[circle,draw=black,thick] (41) at (0,0) {4};
\path[->, thick] (41)edge[] node[]{} (11) (11)edge[] node[]{} (21) (21)edge[] node[]{}(31);
\end{tikzpicture}
         \caption{$t=t_1$}
         \label{fig:small_t1}
     \end{subfigure}
     \hfill
     \begin{subfigure}[b]{0.3\textwidth}
         \centering
\begin{tikzpicture}[]
    \node[circle,draw=black,thick] (12) at (0,1.5) {1};
    \node[circle,draw=black,thick] (22) at (1.5,1.5) {2};
    \node[circle,draw=black,thick] (32) at (1.5,0) {3};
    \node[circle,draw=black,thick] (42) at (0,0) {4};
\path[->, thick] (32)edge[] node[]{} (22) (32)edge[] node[]{} (42);
\end{tikzpicture}
         \caption{$t = t_2$}
         \label{fig:small_t2}
     \end{subfigure}
     \hfill
     \begin{subfigure}[b]{0.3\textwidth}
         \centering
\begin{tikzpicture}[]
    \node[circle,draw=black,thick] (13) at (0,1.5) {1};
    \node[circle,draw=black,thick] (23) at (1.5,1.5) {2};
    \node[circle,draw=black,thick] (33) at (1.5,0) {3};
    \node[circle,draw=black,thick] (43) at (0,0) {4};
\path[->, thick] (43)edge[bend right] node[]{} (13) (13)edge[bend right] node[]{} (43); 
\end{tikzpicture}
         \caption{$t=t_3$}
         \label{fig:small_t3}
     \end{subfigure}
        \caption{Small temporal network ${\mathcal{G}} = (G^{[1]}, G^{[2]},G^{[3]})$.}
        \label{fig:small_graph}
\end{figure}

\subsection{Multilayer interpretation of $M$}\label{sec:multilayer} 
The motivation for the definition of the block matrices above is that we can view walks across multiple time-frames in the following way. 
Recall that given any graph $G$, one can construct the line graph associated with $G$ whose adjacency matrix is the matrix $W$. 
The nodes of the line graph are the edges of $G$ and any two nodes are connected if the concatenation of the corresponding edges in $G$ forms a walk of length two across the original graph. Katz centrality of a node can then be computed by counting walks on the line graph, summing them using appropriate weights, and then using the source and target matrices to ``translate" the  obtained result from the edge space (represented by the line graph) to the node space (represented by $G$). 
The key result that allows this approach is~\eqref{eq:Ar} or, 
equivalently, the fact that a walk of length $r+1$ in the line graph corresponds to a walk of length $r$ in $G$.

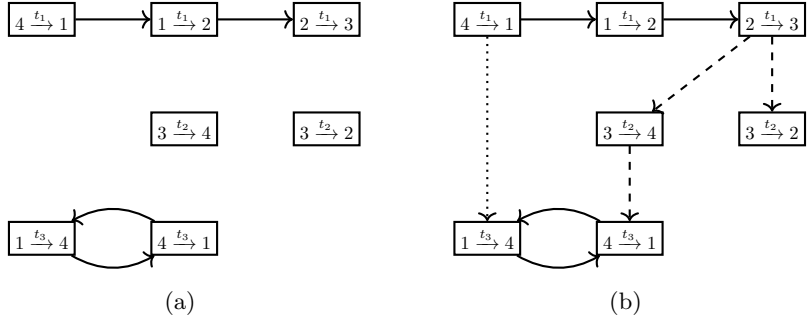
\begin{figure}[]
\centering
\subfloat[]{\label{fig:a}\begin{tikzpicture}[scale=.9]
    \node[scale=.7,rectangle,draw=black,thick] (1) {$ 4 \xrightarrow[]{t_1} 1 $};
    \node[scale=.7,rectangle,draw=black,thick, right=of 1] (2) {$ 1 \xrightarrow[]{t_1} 2 $};
    \node[scale=.7,rectangle,draw=black,thick, right=of 2] (3) {$ 2 \xrightarrow[]{t_1} 3 $};
    \node[scale=.7,rectangle,draw=black,thick, below=of 2] (4) {$ 3 \xrightarrow[]{t_2} 4 $};
    \node[scale=.7,rectangle,draw=black,thick, below=of 3] (5) {$ 3 \xrightarrow[]{t_2} 2 $};
    \node[scale=.7,rectangle,draw=black,thick, below=of 4] (7) {$ 4 \xrightarrow[]{t_3} 1 $};
    \node[scale=.7,rectangle,draw=black,thick, left=of 7] (6) {$ 1 \xrightarrow[]{t_3} 4 $};
\draw[->, thick] (1) to (2);
\draw[->, thick] (2) to (3);
\draw[->, bend right, thick] (6) to (7);
\draw[->, bend right, thick] (7) to (6);
\end{tikzpicture}}
\hspace{1cm}
\subfloat[]{\label{fig:b}\begin{tikzpicture}[scale=.9]
    \node[scale=.7,rectangle,draw=black,thick] (1) {$ 4 \xrightarrow[]{t_1} 1 $};
    \node[scale=.7,rectangle,draw=black,thick, right=of 1] (2) {$ 1 \xrightarrow[]{t_1} 2 $};
    \node[scale=.7,rectangle,draw=black,thick, right=of 2] (3) {$ 2 \xrightarrow[]{t_1} 3 $};
    \node[scale=.7,rectangle,draw=black,thick, below=of 2] (4) {$ 3 \xrightarrow[]{t_2} 4 $};
    \node[scale=.7,rectangle,draw=black,thick, below=of 3] (5) {$ 3 \xrightarrow[]{t_2} 2 $};
    \node[scale=.7,rectangle,draw=black,thick, below=of 4] (7) {$ 4 \xrightarrow[]{t_3} 1 $};
    \node[scale=.7,rectangle,draw=black,thick, left=of 7] (6) {$ 1 \xrightarrow[]{t_3} 4 $};
\draw[->, thick] (1) to (2);
\draw[->, thick] (2) to (3);
\draw[->, bend right, thick] (6) to (7);
\draw[->, bend right, thick] (7) to (6);
\draw[->, dashed, thick] (3) to (4);
\draw[->, dashed, thick] (3) to (5);
\draw[->, dashed, thick] (4) to (7);
\draw[->,dotted, thick] (1) to (6);
\end{tikzpicture}}
\caption{(a) Partial line graph construction associated to ${\mathcal{G}}$ in Figure \ref{fig:small_graph}, as described in subsection~\ref{sec:multilayer}. (b) Multilayer graph underlying the matrix $M$ in Definition~\ref{def:Mnolla} associated with the graph in Figure~\ref{fig:small_graph}}
\label{fig:multi_small_graph}
\end{figure}

To mirror this construction in the case of temporal networks, one begins by constructing the line graphs of each $G^{[\tau]}$ appearing in ${\mathcal{G}}$. 
The new graph representation of ${\mathcal{G}}$ in the edge space will now contain $m=\sum_\tau m_\tau$ nodes, labelled as $i\xrightarrow[]{t_\tau}j$ and hence identified by the source and target of each edge and by the time stamp at which the connection occurs. 
Two nodes in this new graph are connected by a directed edge if their concatenation forms a walk of length two in the graph corresponding to the time stamp $t=t_\tau$, i.e., if they form a walk of length two in $G^{[\tau]}$.  
Consider for example the graph in Figure~\ref{fig:small_graph}. The construction built so far is presented in the left-hand panel of 
Figure~\ref{fig:multi_small_graph}. Clearly, this construction is incomplete, as the walks that we are able to represent, and hence count,  at this stage are only those that take place entirely in one of the $G^{[\tau]}$, for some $\tau = 1,2,\ldots, N$. 
Walks of length two in temporal networks may also take place {\it across} time stamps, and these are not accounted for by the construction of the line graphs $W^{[\tau]}$. 
For example, the temporal walk $4 \xrightarrow[]{t_1} 1 \xrightarrow[]{t_3} 4$ in Figure~\ref{fig:small_graph} is not accounted for in the edge space representation of Figure~\ref{fig:a}.

 The graph built so far, whose adjacency matrix is block-diagonal with $W^{[\tau]}$, $\tau = 1,2,\ldots, N$ as diagonal blocks, thus needs to be adapted to allow for walks that happen across time stamps: directed edges are added to connect nodes of the form $i\xrightarrow[]{t_{\tau_1}} j$ to edges of the form $j\xrightarrow[]{t_{\tau_2} }k$ if $\tau_1<\tau_2$.
  The resulting graph is now a multilayer network which encodes all possible temporal walks of length two taking place in ${\mathcal{G}}$. 
 For example, in the (partial) temporal line-graph construction presented in Figure~\ref{fig:a} associated with the graph in Figure~\ref{fig:small_graph} one needs to add the dashed and dotted edges pictured in Figure~\ref{fig:b}. 
 Here, dashed lines represent connections that happen across two consecutive time stamps, while dotted lines represent walks of length two that start at time $t_1$ and finish at time $t_3$.

If we now turn our attention to the adjacency matrix of this multilayer, it is easy to see that this is indeed the matrix $M$ that we presented in Definition~\ref{def:Mnolla}.

\begin{remark}\label{rem:Mr}
As in the case of static networks, a temporal walk of length $r$ in the graph associated to $M$ corresponds to a temporal walk of length $r+1$ in ${\mathcal{G}}$.
\end{remark}
Because of this interpretation of $M$ as the adjacency matrix of a network associated with ${\mathcal{G}}$, one can expect to be able to compute dynamic walk-based centrality measures using $M$ by counting walks in its associated multilayer graph, summing these walks with appropriate weights, and then projecting the result to the node space using the global source and target matrices.  This was achieved in \eqref{eq:ysym_N} for Katz centrality, and will be generalized to other matrix 
functions in the next section.

\section{$f$-total communicability and $f$-subgraph centrality on time-evolving networks}\label{sec:measures_temporal}
In this section, we consider generalization of the proposed block-matrix approach to walk-based centralities with weights induced by some suitable analytic functions $f$ via their Maclaurin series.

Suppose that $\mathcal{G}$ is a temporal network and let $M=M^{[1,\dots, N]}$ be its global temporal transition matrix. Suppose moreover that, for all $|z|<R$, the series  $f(z) = \sum_{r=0}^\infty c_r z^r$ converges. 
Extending our treatment of the $N=2$ case in section~\ref{sec:dyncom}, we define the temporal $f$-subgraph centrality of node $i$ as 
\[
    \xsym(\alpha)_i = \left(c_0 I_n + \alpha \mathcal{L}^T \sum_{r=1}^\infty c_r (\alpha M)^{r-1} \mathcal{R}\right)_{ii} 
\] 
and the temporal $f$-total communicability of node $i$ as 
\[
\ysym(\alpha)_i = \left(c_0 I_n + \alpha \mathcal{L}^T \sum_{r=1}^\infty c_r (\alpha M)^{r-1} \mathcal{R}\bone_m\right)_{i}.
\]

The following result now follows immediately from the definition of $M$. 
\begin{prop}\label{prop:generalf}
If the power series $f(z) = \sum_{r=0}^\infty c_r z^r$ converges with radius of convergence $R$, then the series $\sum_{r=1}^\infty c_r (\alpha M)^{r-1}$ also converges for all $\alpha\in [0, R/\rho)$, where $\rho = \max_\tau\{\rho(A^{[\tau]})\}$. Moreover,  
\[
\sum_{r=1}^\infty c_r (\alpha M)^{r-1} = \partial f(\alpha M), 
\]
where  $\partial$ is the functional operator \begin{equation}\label{eq:partial}
\partial f(z) := \sum_{r=0}^\infty c_{r+1}z^r = 
\begin{cases}
\frac{f(z)-f(0)}{z} & \text{ if }  z \neq 0\\
f'(0) & \text{ if } z=0.
\end{cases}
\end{equation}
\end{prop}

Using this result we can now formally define temporal $f$-subgraph centrality and temporal $f$-total communicability for any analytic function $f$ as follows.
\begin{definition}\label{def:fcentralities}
Let $\mathcal{G}$ be a temporal network and let $M=M^{[1,\dots, N]}$ be its global temporal transition matrix. Let $f$ be a function which is analytic in a neighborhood of zero with Maclaurin series $f(z) = \sum_{r=0}^\infty c_r z^r$, such that $c_r \geq 0$ for all $r$, and such that $f$ is defined on the spectrum of $M$. Moreover, let $\partial$ be the functional operator defined in \eqref{eq:partial}. {\rm Temporal $f$-subgraph centrality} of node $i$ is defined as
\begin{equation}
    \textbf{x}(\alpha)_i 
= \left(c_0 I +  \alpha\mathcal{L}^T  \partial f(\alpha M) \mathcal{R}\right)_{ii}.
\label{eq:f1}
\end{equation}
Additionally, {\rm temporal $f$-total communicability} of node $i$ is defined as the $i$th entry of the following vector:
\begin{equation}
\textbf{y}(\alpha)
= c_0 \bone_n +  \alpha\mathcal{L}^T \partial f(\alpha M) \bone_m.
\label{eq:f2}
\end{equation}
\end{definition}

These formulae count the number of walks across the temporal network in the required manner: the factor $c_r \alpha^r$ is applied to 
all walks of length $r$. This follows from the one-to-one correspondence between walks of length $r$, weighted as $c_r \alpha^r$ , in the original time-evolving graph and walks of length $r-1$, weighted as $c_r \alpha^{r-1}$, in the graph associated with $M$; see Remark~\ref{rem:Mr}.
Indeed, the role of the operator $\partial$ is to take into account that a walk of length $r-1$ in the multilayer network corresponds to a walk of length $r$ in the original temporal graph, and thus to weight such walks as $c_r \alpha^{r-1}$ (as opposed to $c_{r-1} \alpha^{r-1}$, which arises when $f$ is applied instead of $\partial f$). To recover the 
required weighting, we must then multiply everything by $\alpha$. Finally, we add\footnote{Of course, adding such a term does not change the ranking associated with the corresponding $f$-centrality,  and therefore it is not necessary in a practical implementation of an algorithm to compute it. We incorporate it to obtain the correct combinatorial expressions.} a term $c_0 I$ since walks of length zero in the original time-evolving graph do not have any correspondence in the graph associated to $M$. Observe that in the special case of Katz centrality $\partial f(z)=f(z)=(1-z)^{-1}$, and therefore the importance of applying the operator $\partial$ is not manifest.

\begin{example}\label{ex:expono}
Let us consider $f(z)=e^z$, so that \[ \partial f(z) = \psi_1(z) = \sum_{r=0}^\infty \frac{z^r}{(r+1)!},\]
in the case of 
three 
time stamps. Then, the exponential subgraph centrality for all $\beta>0$ can be computed as the diagonal entries of
\begin{equation}\label{eq:genfunexp3}
I + \beta \mathcal{L}^T \psi_1 \left( \beta \begin{pmatrix} W^{[1]} & W^{[1,2]} & W^{[1,3]}\\
0&W^{[2]}&W^{[2,3]}\\
0&0&W^{[3]}  \end{pmatrix} \right)  \mathcal{R} 
\end{equation}  
while the total communicability vector is obtained by multiplying \eqref{eq:genfunexp3} by $\mathbf{1}$.
We emphasize that the generating function \eqref{eq:genfunexp3} is generally \emph{not} equal to $e^{\beta A^{[1]}} e^{\beta A^{[2]}}e^{\beta A^{[3]}}$. 

Indeed, consider for example the temporal graph ${\mathcal{G}}$ represented by the following adjacency matrices
\[
A^{[1]}= \begin{bmatrix}
0 & 1 & 0 & 0 \\ 
0 & 0 & 0 & 0 \\
0 & 0 & 0 & 0 \\
0 & 0 & 0 & 0 
\end{bmatrix}, \, 
A^{[2]}= \begin{bmatrix}
0 & 0 & 0 & 0 \\ 
0 & 0 & 1 & 0 \\
0 & 0 & 0 & 0 \\
0 & 0 & 0 & 0 
\end{bmatrix}, \text{ and } 
A^{[3]}= \begin{bmatrix}
0 & 0 & 0 & 0 \\ 
0 & 0 & 0 & 0 \\
0 & 0 & 0 & 1 \\
0 & 0 & 0 & 0 
\end{bmatrix}.
\]
Clearly, all the admissible walks in this temporal graph are 
\begin{itemize}
    \item walks of length one, weighted by $\beta$:
    \[1 \xrightarrow[]{t_1} 2, \; 2 \xrightarrow[]{t_2} 3,\; \text{ and } \; 3 \xrightarrow[]{t_3} 4;\]
    \item walks of length two, weighted by $\frac{\beta^2}{2}$:
    \[1 \xrightarrow[]{t_1} 2 \xrightarrow[]{t_2} 3\; \text{and }\; 2 \xrightarrow[]{t_2} 3 \xrightarrow[]{t_3} 4;\]
    and finally, 
    \item one walk of length three, weighted by a factor $\frac{\beta^3}{3!}$
\[1 \xrightarrow[]{t_1} 2 \xrightarrow[]{t_2} 3 \xrightarrow[]{t_3} 4.\]
\end{itemize}
Since these are the only admissible walks, together with walks of length zero, weighted by $1$, 
we expect the entries of the $4\times 4$ matrix~\eqref{eq:genfunexp3} to encode these weights in the appropriate entries; 
for example, in the entry $(1,3)$ we expect to read the weight $\frac{\beta^2}{2}$ since the only walk in the 
network originating at node $1$ and terminating at node $3$ is of length two. 

It is easy to verify that, up to permutation similarity, 
\[
M = \begin{bmatrix}
0 & 1 & 0 \\ 0 & 0 & 1 \\ 0 & 0 & 0
\end{bmatrix}
\]
which is nilpotent, and hence 
\[\psi_1(\beta M) = I_3 + \frac{\beta}{2} M + \frac{\beta^2}{3!}M^2 = 
\begin{bmatrix}
1 & \frac{\beta}{2} & \frac{\beta^2}{3!} \\
0 & 1 & \frac{\beta}{2} \\
0 & 0 & 1
\end{bmatrix}.\]
With the labelling of edges induced by $M$ above it holds that
\[
{\mathcal{L}} = 
\begin{bmatrix}
1 & 0 & 0 & 0 \\
0 & 1 & 0 & 0 \\ 
0 & 0 & 1 & 0
\end{bmatrix}\quad \text{ and }\quad
{\mathcal{R}} =
\begin{bmatrix}
0 & 1 & 0 & 0 \\
0 & 0 & 1 & 0 \\ 
0 & 0 & 0 & 1
\end{bmatrix}.
\]
It immediately follows that the matrix in \eqref{eq:genfunexp3} is 
\[
I_4 + \beta {\mathcal{L}}^T\psi_1(\beta M){\mathcal{R}} = 
\begin{bmatrix}
1 & \beta & \frac{\beta^2}{2} & \frac{\beta^3}{3!} \\
0 & 1 & \beta & \frac{\beta^2}{2} \\
0 & 0 & 1 & \beta \\
0 & 0 & 0 & 1
\end{bmatrix},
\]
which encodes in its entries exactly the required weights. 
On the other hand, multiplying the exponentials of the adjacency matrices gives 
\[
e^{\beta A^{[1]}}e^{\beta A^{[2]}}e^{\beta A^{[3}} = 
\begin{bmatrix}
1 & \beta & \beta^2 &\beta^3 \\
0 & 1 & \beta & \beta^2 \\
0 & 0 & 1 & \beta \\
0 & 0 & 0 & 1
\end{bmatrix},
\]
where we have used $(A^{[\tau]})^r = 0$ for all $r\geq 2$ and  $A^{[1]}A^{[3]} = 0$. 
This matrix fails 
to correctly weight walks of length two and three. 
\end{example}
\smallskip

In the next section we move on to the problem 
of generalizing the results obtained so far to the 
non-backtracking setting, where backtracking is 
disallowed either within time stamps or across 
time stamps. 
The construction presented mirrors the one described 
above and will build on modified versions of the global 
temporal transition matrix $M$.

\section{Non-backtracking walks}\label{sec:nbt} 

Non-backtracking walks, as defined in
Definition~\ref{def:nbt}, 
may be studied using the 
{\it Hashimoto matrix}; see~\cite{Hash90}.
We also note that the 
Hashimoto matrix plays a key role in related studies of 
\emph{random} NBT walks \cite{FTT21}.
\begin{definition}
The {\rm Hashimoto matrix} associated with a graph $G=(V,E)$ with nodes $V = \{1,2,\ldots, n\}$ and $m$ edges is the matrix $B\in\R^{m\times m}$ defined as
\[
B_{i\to j,k\to \ell} = \delta_{jk}(1-\delta_{i\ell})
\] 
for all $i\to j,k\to \ell\in E$. 
\end{definition}
It is easily verified that 
\[
B = W - W\circ W^T
\]
where $\circ$ denotes the Schur (or entrywise) product and $W$ is the adjacency matrix of the 
line graph of $G$.
In words, rows and columns of $B$ correspond to edges in $G$, and 
the nonzero entries record pairs of edges that form a 
NBT walk of length two in $G$.

Taking powers of this matrix and then projecting back to the node space using the source and target 
matrices allows us to count NBT walks taking place in $G$. Specifically, 
\begin{equation}\label{eq:Br}
(L^TB^{r-1}R)_{ij} = (p_r(A))_{ij}
\end{equation}
is the number of NBT walks of length $r$ originating at node $i$ and ending at node $j$; see~\cite{Beyond}.
Here and in the following we denote by $p_r(A)\in\Rnn$ the matrix whose entries record the number 
of NBT walks of length $r$ between any two nodes in the network with adjacency matrix $A$. 
These matrices satisfy a four-term recurrence that allows the computation of the number 
of NBT walks of any length in the network. 
We refer the interested reader to~\cite{AGHN17a,GHN18} and references therein.
 
\subsection{Non-backtracking walk-based centrality measure}
A NBT version of Katz centrality was proposed and studied in~\cite{AGHN17a,GHN18}. 
Generalization to other $f$-centrality measures was presented in~\cite{AGHN17b}. 
Later, in~\cite{Beyond} it was shown how NBT Katz centrality, as well as other NBT versions of $f$-subgraph centrality and $f$-total communicability, could also be derived via the Hashimoto matrix.

We may  
define the NBT Katz centrality by changing 
(\ref{eq:Katz}) to 
\[
\widehat{\yvec}(\alpha) =
\sum_{r=0}^\infty \alpha^r  p_r(A)\bone.
\]
  
As shown in \cite{Beyond}, we may also obtain $\widehat{\yvec}(\alpha)$ 
by computing Katz centrality for the graph of the Hashimoto matrix and then projecting via the source and target matrices; that is, 
\[
    \widehat{\yvec}(\alpha) = \left(I + \alpha L^T(I-\alpha B)^{-1}R\right) \bone, 
\]
for $\alpha<1/\rho(B)$. 
Similarly, NBT $f$-total communicability is obtained by computing $\partial f$-total communicability for the Hashimoto matrix and then projecting, where $\partial$ is the operator defined in \eqref{eq:partial}, 
that is, 
\begin{equation*}
    \widehat{\yvec}(\alpha) = \left(f(0)I + \alpha L^T \partial f (\alpha B) R\right) \boldsymbol{1}.
\end{equation*}
NBT versions of $f$-subgraph centrality can also be defined analogously using the diagonal elements of $f(0)I + \alpha L^T \partial f (\alpha B) R$.

\subsection{Temporal NBT centrality measures}

In the remainder of this section we 
carry through concepts and results
from sections~\ref{sec:dyncom} and \ref{sec:measures_temporal}
to the 
NBT setting.
This produces, for the first time,
definitions and computable 
expressions for temporal NBT centrality measures. 
We first notice that there is not just one type of backtracking for temporal networks; indeed, there are three:
\begin{itemize}
\item Backtracking happens within a certain time stamp; we will refer to this as {\it backtracking in space}.
\item Backtracking happens across time stamps; we will refer to this as {\it backtracking in time}.
\item Backtracking happens both within a time stamp and across time stamps (not necessarily in that order); we will refer to this as {\it backtracking in time and space}.
\end{itemize}
To explain these different types, consider the temporal graph ${\mathcal{G}}$ in Figure~\ref{fig:small_graph}. 
The walks $1 \xrightarrow[]{t_3} 4 \xrightarrow[]{t_3} 1$ and $2 \xrightarrow[]{t_1} 3 \xrightarrow[]{t_2} 2$ showcase two different types of backtracking: 
the first is backtracking in space and the second is backtracking in time; additionally, the walk $4 \xrightarrow[]{t_1} 1 \xrightarrow[]{t_3} 4 $ provides an example of a walk that backtracks in both space \textit{and} time.

 \subsubsection{Non-backtracking Katz centrality for time-evolving networks: the case of $2$ time-frames}

As in section~\ref{sec:measures_temporal}, we begin by considering the case of $\mathcal{G} = (G^{[1]},G^{[2]})$ at times $t_1 < t_2$.  
In the following, we make use of the notation $B^{[1,2]} := W^{[1,2]} - W^{[1,2]}\circ {W^{[2,1]}}^T$. 
\begin{remark}
Recall that the matrix $W^{[1,2]} = R^{[1]}(L^{[2]})^T\in\R^{m_1\times m_2}$ is a possibly rectangular matrix that encodes in its entries whether two edges, one appearing at time $t=t_1$ and a second appearing at time $t=t_2$, can be concatenated; see Definition~\ref{def:Wt1t2} and Proposition~\ref{prop:Wt1t2}. 
Analogously, the matrix 
$B^{[1,2]}\in\R^{m_1\times m_2}$ records whether a walk of length two that occurs across time stamps 
$t_1$ and $t_2$ is non-backtracking, since we have 
\[
(B^{[1,2]})_{i\to j,k\to \ell} = \delta_{jk}(1-\delta_{i\ell})
\]
for all $i\to j\in E^{[1]}$ and for all $k\to\ell\in E^{[2]}$. 
We will formalize this concept in Definition~\ref{def:Bt1t2} and Proposition~\ref{prop:Bt1t2} below.
\end{remark}

With consideration to the previous interpretation of the terms appearing in~\eqref{eq:goodexample}, if one wants to forbid backtracking, in full or in part, it is appropriate to generalize the expression~\eqref{eq:goodexample}  in three possible ways, depending on the type of backtracking behavior that we want to forbid.
\begin{itemize}
\item Backtracking in space forbidden:  we replace the $W^{[\tau]}$ with the Hashimoto matrix $B^{[\tau]}$ for $\tau=1,2$, giving  
\begin{eqnarray}
I_n   &+& \alpha \left(  \sum_{\tau=1}^2 (L^{[\tau]})^T (I_m - \alpha B^{[\tau]})^{-1} R^{[\tau]}      \right) \nonumber\\
&&
+ \alpha^2 (L^{[1]})^T (I_m-\alpha B^{[1]})^{-1} W^{[1,2]} (I_m- \alpha B^{[2]})^{-1} R^{[2]}. \label{eq:one}
\end{eqnarray} 
\item Backtracking in time forbidden: we replace $W^{[1,2]}$ with $B^{[1,2]}$, giving 
\begin{eqnarray}
I_n   &+& \alpha \left(  \sum_{\tau=1}^2 (L^{[\tau]})^T (I_m - \alpha W^{[\tau]})^{-1} R^{[\tau]}      \right) \nonumber \\
&&
+ \alpha^2 (L^{[1]})^T (I_m-\alpha W^{[1]})^{-1} B^{[1,2]} (I_m- \alpha W^{[2]})^{-1} R^{[2]}. 
\label{eq:two}
\end{eqnarray}  
\item Backtracking in time and space both forbidden: we perform both replacements, giving 
\begin{eqnarray}
I_n   &+& \alpha \left(  \sum_{\tau=1}^2 (L^{[\tau]})^T (I_m - \alpha B^{[\tau]})^{-1} R^{[\tau]}      \right) \nonumber \\
&& + \alpha^2 (L^{[1]})^T (I_m-\alpha B^{[1]})^{-1} B^{[1,2]} (I_m- \alpha B^{[2]})^{-1} R^{[2]}.  \label{eq:three}
\end{eqnarray} 
\end{itemize}

In the case of \eqref{eq:one}, we observe that for $\tau = 1,2$   
\[ \alpha(L^{[\tau]})^T (I_m - \alpha B^{[\tau]})^{-1} R^{[\tau]} = \alpha\sum_{r=0}^\infty \alpha^r p_{r+1}(A^{[\tau]}), \] 
where we have used \eqref{eq:Br}.

Moreover, for small enough values of $\alpha$,
\[ \sum_{r=0}^\infty \alpha^r p_r(A^{[\tau]}) = (1-\alpha^2) [I - \alpha A^{[\tau]} + \alpha^2 (D^{[\tau]} - I) + \alpha^3 (A^{[\tau]}-S^{[\tau]})]^{-1}, \] where
\begin{equation}\label{eq:dands}
D^{[\tau]} = \diag \diag (A^{[\tau]})^2, \qquad S^{[\tau]} = A^{[\tau]} \circ (A^{[\tau]})^T
\end{equation}
 (both formulae are proved in the existing literature \cite{AGHN17a,Beyond}). 
 Hence, all the terms appearing in \eqref{eq:one} can be reformulated 
 in terms of matrices of order $n$,
 thus obtaining a \emph{non-backtracking (in space only) dynamic communicability matrix}:
\begin{equation*}\label{eq:nbtisodcm}
 \mathcal{Q} = (1 - \alpha^2)^2 \prod_{\tau=1}^2 [I - \alpha A^{[\tau]} + \alpha^2 (D^{[\tau]} - I) + \alpha^3 (A^{[\tau]}-S^{[\tau]})]^{-1}. 
\end{equation*}  

When possible, it is clearly preferable to work at node-level, since for large values of $n$ it may easily happen that $n \ll m$. We observe that this may be true even if, individually, each graph has $O(n)$ edges, since in total there may be up to $O(Nn)$ edges. However, a direct node-level formulation appears to be infeasible when backtracking is also forbidden in time, due to the fact that a memory of the last traveled edge in a given time frame is necessary. Thus, below we will focus on developing further the edge level formulations \eqref{eq:one}--\eqref{eq:three}.

Our goal is to generalize these formulae to the case of $N>2$ time stamps. 
When doing so, one issue that arises is that the above formulae for the dynamic communicability matrices become unwieldy. 
The solution we propose builds on the fact that \eqref{eq:one}--\eqref{eq:three} can be rewritten using block matrices, just like we did with~\eqref{eq:goodexample}. 
Indeed, depending on what type of backtracking is forbidden, one can generalize \eqref{eq:Mzero} as follows:
\begin{equation}\label{eq:M2general}
M^{[1,2]} = \begin{bmatrix} C^{[1]} & C^{[1,2]}\\
0 & C^{[2]}
\end{bmatrix},
\end{equation}
where
\begin{itemize}
    \item $C^{[\tau]} = B^{[\tau]}$ for $\tau = 1,2$ and $C^{[1,2]} = W^{[1,2]}$, if one wants to recover \eqref{eq:one}, i.e., avoid backtracking in space but not in time; 
    \item $C^{[\tau]} = W^{[\tau]}$ for $\tau = 1,2$ and $C^{[1,2]} = B^{[1,2]}$, if one wants to recover \eqref{eq:two}, i.e., avoid backtracking in time but not in space; and
    \item $C^{[\tau]} = B^{[\tau]}$ for $\tau = 1,2$ and $C^{[1,2]} = B^{[1,2]}$, if one wants to recover \eqref{eq:three}, i.e., avoid backtracking in space and time.
\end{itemize}
The following result then holds.
\begin{prop}\label{prop:properblocks}
Suppose $M^{[1,2]}$ is defined as in~\eqref{eq:M2general}. Then, within the radius of convergence,
  \[ (I_{m} - \alpha {M^{[1,2]}})^{-1}  = \begin{pmatrix}(I_m-\alpha C^{[1]})^{-1} & (I_m-\alpha C^{[1]})^{-1}\alpha C^{[1,2]}(I_m-\alpha C^{[2]})^{-1} \\ 0 & (I_m-\alpha C^{[2]})^{-1}\end{pmatrix}.\]
\end{prop}
Our proof of this result mirrors the proof of Proposition ~\ref{prop:properblocksW} and is therefore omitted.

By using the above proposition and the fact that 
\[(I_m-\alpha M^{[1,2]})^{-1} = I_{m} + \sum_{r = 1}^\infty \alpha^r {M^{[1,2]}}^r,\] 
we can easily recover~\eqref{eq:one}--\eqref{eq:three} by projecting back to the node-level, multiplying by $\alpha$ to account for the fact that a walk on the edge level is one step longer than its equivalent on the node level, and adding back the identity matrix $I_n$; such a projection is achieved by the global source and target matrices associated with ${\mathcal{G}}$: 
\[
 I_n + \alpha {\mathcal{L}^{[1,2]}}^T (I_m-\alpha M^{[1,2]})^{-1} \mathcal{R}^{[1,2]}.\]

\smallskip

We now move to the goal of extending this line-graph interpretation of Katz centrality for evolving networks to possibly more than two graphs. In doing so, we will allow for backtracking to be forbidden (in time or space or both).
We adopt the same strategy used in sections~\ref{sec:dyncom} and \ref{sec:measures_temporal}.

\subsubsection{Non-backtracking Katz centrality for time-evolving networks: more than $2$ time-frames}
We begin by extending Definition~\ref{def:Wt1t2} and Definition~\ref{def:Mnolla} to the NBT framework. The latter 
may be reformulated in different ways, depending on the type of backtracking that one wants to avoid. 
\begin{definition}\label{def:Bt1t2}
Let ${\mathcal{G}} = (G^{[1]},G^{[2]},\ldots,G^{[N]})$ be a temporal network with $N$ time stamps. Let moreover $G^{[\tau]} = ({\mathcal{V}}, E^{[\tau]})$ for all $\tau = 1,2,\ldots, N$.  
The {\rm non-backtracking transition matrix} $B^{[\tau_1,\tau_2]}$ representing walks of length two that traverse the first edge at time $t_{\tau_1}$ and the second at time $t_{\tau_2}$, $\tau_1 < \tau_2$, is defined entry-wise as follows:
\[
(B^{[\tau_1,\tau_2]})_{i \rightarrow j, k \rightarrow \ell} = \delta_{jk}(1-\delta_{i\ell})
\]
for all $i\to j\in E^{[\tau_1]}$ and $k\to\ell\in E^{[\tau_2]}$. 
\end{definition}

These matrices, which are generally non-square, record information about NBT walks of length two where the first edge is traversed at time $\tau_1$ and the second is traversed at time $\tau_2$ (after staying idle for some time, in the case where $\tau_2 > \tau_1+1$). 
The next result extends Proposition~\ref{prop:Wt1t2}.
\begin{prop}\label{prop:Bt1t2}
Let ${\mathcal{G}} = (G^{[1]},G^{[2]},\ldots,G^{[N]})$ be a temporal network. 
Let moreover $L^{[\tau]}$ and $R^{[\tau]}$ be the source and target matrices associated with graph $G^{[\tau]}$. 
Then for all $1\leq\tau_1<\tau_2\leq N$ we have  
\[
    B^{[\tau_1,\tau_2]} = W^{[\tau_1,\tau_2]} - W^{[\tau_1,\tau_2]}\circ {W^{[\tau_2,\tau_1]}}^T= R^{[\tau_1]}(L^{[\tau_2]})^T - (R^{[\tau_1]}L^{[\tau_2]^T}) \circ (L^{[\tau_1]}R^{[\tau_2]^T}).
\]
\end{prop}
\begin{proof}
We proceed entrywise. Let $i\to j\in E^{[\tau_1]}$ and $k\to\ell\in E^{[\tau_2]}$ 
be two edges in ${\mathcal{G}}$. Then
\[
(R^{[\tau_1]}L^{[\tau_2]^T})_{i\to j,k\to \ell} = 
\sum_{h = 1}^n (R^{[\tau_1]})_{i\to j,h} (L^{[\tau_2]})_{k\to \ell, h}= \delta_{jk}.
\]
Similarly,
\[
(L^{[\tau_1]}R^{[\tau_2]^T})_{i\to j,k\to\ell} = 
\sum_{h = 1}^n (L^{[\tau_1]})_{i\to j,h} (R^{[\tau_2]})_{k\to \ell, h}=  \delta_{i\ell}.
\]
The conclusion then follows from Definition~\ref{def:Bt1t2} and Proposition~\ref{prop:Wt1t2}.

\end{proof}

An advantage of the block-matrix approach is that it is readily expanded to deal with $N$ time-frames. This is done by defining $M^{[1,\dots,N]}$ as one of the following block upper triangular block matrices.
\begin{definition}\label{def:Mkolme}
Let ${\mathcal{G}} = (G^{[1]},G^{[2]},\ldots,G^{[N]})$ be a temporal network with $N$ time stamps. 
The {\rm NBT global temporal transition matrix} associated with ${\mathcal{G}}$ is the $m\times m$ block matrix 
\begin{equation}\label{eq:Mkolme}
M = M^{[1,\dots,N]} = 
\begin{bmatrix}   
C^{[1]} & C^{[1,2]} & C^{[1,3]} & \dots & C^{[1,N]} \\ 
0 & C^{[2]} & C^{[2,3]} & \dots & C^{[2,N]} \\
\vdots & & \ddots &  & \vdots \\
\vdots & & & \ddots & \vdots\\
0 & \dots & \dots & 0 & C^{[N]}
\end{bmatrix}, 
\end{equation}
where 
\begin{itemize}
\item[(i)] $C^{[\tau_1]} = B^{[\tau_1]}$ and $C^{[\tau_1,\tau_2]} = W^{[\tau_1,\tau_2]}$ for all $\tau_1,\tau_2=1,2,\ldots,N$ ($\tau_1< \tau_2$) if backtracking in space is forbidden;
\item[(ii)] $C^{[\tau_1]} = W^{[\tau_1]}$ and $C^{[\tau_1,\tau_2]} = B^{[\tau_1,\tau_2]}$ for all $\tau_1,\tau_2=1,2,\ldots,N$ ($\tau_1< \tau_2$) if backtracking in time is forbidden; and 
\item[(iii)] $C^{[\tau_1]} = B^{[\tau_1]}$ and $C^{[\tau_1,\tau_2]} = B^{[\tau_1,\tau_2]}$ for all $\tau_1,\tau_2=1,2,\ldots,N$ ($\tau_1< \tau_2$) if backtracking in both time and space is forbidden.
\end{itemize} 
\end{definition}

Such a formulation allows for non-backtracking across multiple graphs to be calculated.

\begin{remark}
In the third case, one can calculate $M$ from the global source and target matrices as $\mathcal{R}\mathcal{L}^T - \mathcal{R}\mathcal{L}^T \circ \mathcal{L}\mathcal{R}^T$ by setting all elements below the block diagonal to zero.
\end{remark}

We can now define Katz centrality on a time-evolving network with NBT global temporal transition matrix $M=M^{[1,\dots, N]}$ equal to any of the three options in Definition~\ref{def:Mkolme} (depending on the type of backtracking allowed)
as
\[ \widehat{\ysym}(\alpha) = \sum_{r=0}^\infty (\alpha M)^r \bone.\]

The following result gives an upper bound for the possible choices of $\alpha>0$.
\begin{prop}\label{prop:radius}
Let ${\mathcal{G}} = (G^{[1]},G^{[2]}, \ldots, G^{[N]})$ be a temporal network with 
adjacency matrices $A^{[1]}, A^{[2]}, \ldots, A^{[N]}$ and let $M = M^{[1,\dots, N]}$ be 
defined as in Definition~\ref{def:Mkolme}. 
Moreover, let $\rho=\max_\tau \rho(A^{[\tau]})$ and $\lambda=\min_\tau \lambda_\tau$,  where $\lambda_\tau$ is the smallest eigenvalue of the matrix polynomial $P^{[\tau]}(z)=I-A^{[\tau]} z + (D^{[\tau]}-I) z^2 + (A^{[\tau]}-S^{[\tau]}) z^3$, for all $\tau = 1,2,\ldots, N$, where $D^{[\tau]}$ and $S^{[\tau]}$ are as in \eqref{eq:dands}.  

Then, the series $\sum_{r=0}^\infty (\alpha M)^r$ converges and moreover
\[
\sum_{r=0}^\infty (\alpha M^{[1,\dots,N]})^r  = (I_m-\alpha M^{[1,\dots, N]})^{-1}
\]
for all $\alpha\in (0,\ell)$,  where
\begin{itemize}
\item $\ell=\lambda$, if $M$ is as in {\it (i)} or {\it (iii)} from Definition~\ref{def:Mkolme}; %
\item $\ell=\rho^{-1}$, if $M$ is as in {\it (ii)} from Definition~\ref{def:Mkolme}. 
\end{itemize}
\end{prop}
\begin{proof}
The rightmost endpoint of the interval of allowed values for $\alpha$ is the radius of convergence of $\sum_{r=0}^\infty M^r$, which is equal to the inverse of the spectral radius of $M$. Since the latter is a block triangular matrix, its spectral radius is equal to the maximum of the spectral radii of the diagonal blocks. It follows by the analysis in \cite{AGHN17a, Beyond, GHN18} that the spectral radius of $B^{[i]}$ is $\lambda^{-1}_i$, while by Flanders' theorem the spectral radius of $W^{[i]}=R^{[i]}(L^{[i]})^T$ is equal to the spectral radius of $A^{[i]}=(L^{[i]})^T R^{[i]}$, which is $\rho_i$.
\end{proof}

This convergence result allows us to thus define NBT Katz centrality for temporal networks as follows.
\begin{definition}\label{def:katznbt}
Let ${\mathcal{G}} = (G^{[1]},G^{[2]},\ldots,G^{[N]})$ be a temporal network with $N$ time stamps.
Let $M$ be defined as in Definition~\ref{def:Mkolme} and let $\alpha\in (0,\ell)$, where $\ell$ is defined according to Proposition~\ref{prop:radius}.
We define {\rm NBT temporal Katz centrality} as the vector 
\[
 \widehat{\ysym}(\alpha) = \sum_{r=0}^\infty (\alpha M)^r \bone = (I_n + \alpha (\mathcal{L}^T (I_m-\alpha M)^{-1} \mathcal{R}) \bone_n.
\]
\end{definition}

\begin{remark}
Note that, since $\mathcal{R}\bone_n=\bone_m$, Katz centrality is also equal to
\begin{equation*}\label{eq:k1}
 \widehat{\ysym}(\alpha) = \bone_n+ \alpha \mathcal{L}^T (I-\alpha M)^{-1} \bone_m.
\end{equation*}
\end{remark}
\begin{remark}
One may also define subgraph centrality indices by considering the 
diagonal entries of the matrix under study:
\[
\widehat{\xsym}(\alpha)_i = 1 + \alpha  (\mathcal{L}^T (I_m-\alpha M)^{-1} \mathcal{R})_{ii}.
\]
\end{remark}

The results presented so far for the matrix resolvent can also be extended to treat other analytic matrix functions. Indeed, the following result is an analogue of Proposition~\ref{prop:generalf} in the non-backtracking framework. 
\begin{prop}\label{prop:fnbt}
If the power series $f(z) =\sum_{r=0}^\infty c_rz^r$ converges with radius of convergence $R$,  then  the  series
$\sum_{r=0}^\infty c_r \alpha^r M^r$ also  converges  for  all $\alpha \in [0,R/\ell)$, with $M$ is as in Definition~\ref{def:Mkolme} and  
$\ell$ is as in Proposition~\ref{prop:radius}.  
Moreover, 
\[\sum_{r=0}^\infty c_r \alpha^r M^r =\partial f(\alpha M),\] where $\partial$ is the functional operator defined in \eqref{eq:partial}.
\end{prop}
Using this result, Definition~\ref{def:katznbt} can straightforwardly be generalized to describe centrality indices that avoid different types of backtracking, according to how the matrix $M$ in \eqref{eq:Mkolme} is selected. 
\begin{definition}
Let $\mathcal{G}$ be a temporal network and let $M=M^{[1,\dots, N]}$ be defined as in Definition~\ref{def:Mkolme} Let $f$ be a function which is analytic in a neighborhood of zero with Maclaurin series $f(z) = \sum_{r=0}^\infty c_r z^r$, such that $c_r \geq 0$ for all $r$, and let $\alpha\in(0,\ell)$, where $\ell$ is as in Proposition~\ref{prop:radius}. 
Moreover, let $\partial$ be the functional operator defined in \eqref{eq:partial}. 
The {\rm NBT temporal $f$-subgraph centrality} of node $i$ is defined as
\begin{equation}
    \widehat{\xsym}(\alpha)_i 
= \left(c_0 I +  \alpha\mathcal{L}^T  \partial f(\alpha M) \mathcal{R}\right)_{ii}.
\label{eq:nbtf1}
\end{equation}
Additionally, the {\rm NBT temporal $f$-total communicability} of node $i$ is defined as the $i$th entry of the following vector:
\begin{equation}
\widehat{\ysym}(\alpha)
= c_0 \bone_n +  \alpha\mathcal{L}^T \partial f(\alpha M) \bone_m.
\label{eq:nbtf2}
\end{equation}
where $\mathcal{L}$ and $\mathcal{R}$ be defined as in Definition~\ref{def:globalLR}. 
\end{definition}

\section{Conclusions}\label{sec:conc}

Our aim in this work was to study combinatorically-motivated centrality measures for temporal networks.
This required us to extend existing work involving
matrix functions applied to (a) the adjacency matrix of a network,
(b) the adjacency matrix of the line graph, and (c) the Hashimoto matrix.
For the Katz, or matrix resolvent, case, it is known 
from \cite{GHPE11} that dynamic 
centrality may be computed at the node level; 
that is, by using $n$ by $n$ matrices, where $n$ is the number of nodes. 
Here, we proved that any other matrix function-based measure
does not have such a straightforward 
representation, but can be expressed in terms of a higher
dimensional edge-level formula.
In the case of non-backtracking temporal centrality,
we identified three types of constraint 
and showed how these may be computed from block matrix 
extensions of the 
Hashimoto construction.

 ------------------------------------------------------------------

\bibliography{netrefs}

\bibliographystyle{amsplain}

\end{document}